\documentstyle[11pt]{article}
\newtheorem {thm}{Theorem}[section]
\newtheorem {prop}[thm]{Proposition}
\newtheorem {lem}[thm]{Lemma}
\newtheorem {cor}[thm]{Corollary}
\newtheorem {defn}[thm]{Definition}

\input{epsf.sty}
\def\Cox{\hfill \Box}
\def\Z{{\bf Z}}
\def\P{{\bf P}}
\def\R{{\bf R}}

\def\0{{\bf 0}}

\def\eps{\varepsilon}

\def\cut{{\rm cut}}
\def\ext{{\rm ext}}
\def\ERB{\texttt{ERB}}
\def\ELT{\texttt{ELT}}
\def\RF{\texttt{RF}}

\begin{document}

\title{Some two-dimensional finite energy percolation
processes\thanks{Research of the first author supported by the Swedish
Research Council and by the G\"oran Gustafsson Foundation for Research in the
Natural Sciences and Medicine. Research of the second author partially
supported by OTKA grant K49398.}}
\author{Olle H\"aggstr\"om\thanks{Dept.\ of Mathematical Sciences,
Chalmers University of Technology, 412 96 G\"oteborg, Sweden,
\texttt{http://www.math.chalmers.se/\~{ }olleh/}.}
\and P\'eter Mester\thanks{Dept.\
of Mathematics, Indiana University,
Bloomington, IN 47405-5701, United States.}}
\maketitle

\begin{abstract}
Some examples of translation invariant site percolation processes
on the $\Z^2$ lattice are constructed, the most far-reaching example being
one that satisfies uniform finite
energy (meaning that the probability that a site is open given the
status of all others is bounded away from $0$ and $1$) and exhibits a.s.\
the coexistence of an infinite open cluster and an infinite closed cluster.
Essentially the same example shows that coexistence is possible
between an infinite open cluster and an infinite closed cluster that are
both robust under i.i.d.\ thinning. 
\end{abstract}

\section{Introduction}  \label{sect:intro}

By a {\bf site percolation} on $\Z^2$, we mean an
$\{0,1\}^{\Z^2}$-valued random object $X$.
Focus in percolation theory is mainly on the connected components
(clusters) of $X$. Two vertices $x,y\in \Z^2$ are said to
communicate if there exists a path
$\{z_1, z_2, \ldots, z_n\}$ from $z_1=x$ to $z_n=y$ with
$X(z_1)=X(z_2)=\cdots = X(z_n)$
(in the definition of path, we require $z_i$ and $z_j$ to be
$L_1$-nearest neighbors for each $i$), and a connected component is
a maximal set of vertices that all communicate with each other.
A connected component is called an {\bf open cluster} or a {\bf closed
cluster} depending on whether its vertices take value $1$ or $0$ in $X$.

Much of percolation theory deals with the i.i.d.\ case
(see, e.g., Grimmett \cite{G}), though various dependent settings
have also received much attention. Here we will abandon
the i.i.d.\ assumption in favor of the weaker but natural assumption
of {\bf translation invariance}, meaning that for any
$n$ and any $x_1, \ldots, x_n \in \Z^2$, the distribution of
$(X(x_1+y), \ldots, X(x_n+y))$ does not depend on $y \in \Z^2$.

Intuitively,
the planar structure of $\Z^2$ makes it difficult for an infinite open
cluster and an infinite closed cluster to coexist. In order
to prove a theorem
to this extent, some further conditions beyond translation invariance
are needed, as the following trivial example shows: assign the
vertices on the $x$-axis independently value $0$ or $1$ with
probability $\frac{1}{2}$ each, and let the values of any other vertex $z$
be dictated by the value of the vertex on the $x$-axis sharing $z$'s
$x$-coordinate. This produces a translation invariant site percolation
with both infinite open clusters and infinite closed clusters
(in fact, infinitely many of each).

In a seminal paper,
Gandolfi, Keane and Russo \cite{GKR}
showed that translation invariance and positive associations together with
some auxiliary ergodicity conditions (later relaxed by Sheffield \cite{Sh})
is enough to rule out such coexistence. In applications of this result,
the hard part has typically been to establish positive associations; see for
instance Chayes \cite{C} and H\"aggstr\"om \cite{H}.
Partly for this reason, several
researchers over the years have asked whether the positive
associations condition can be replaced by the often easier-to-verify
condition of {\bf finite energy}, defined as follows.
\begin{defn}
A site percolation $X$ on $\Z^2$ is said to satisfy {\bf finite energy}
if it admits conditional probabilities such that for all $x\in \Z^2$
and all $\xi \in \{0,1\}^{\Z^2 \setminus \{x\}}$ we have
\[
0 \, < \, \P(X(x)=1 \, | \, X(\Z^2 \setminus \{x\}) = \xi) \, < \, 1 \, .
\]
It is said to satisfy {\bf uniform finite energy}
if for some $\eps>0$
it admits conditional probabilities such that for all $x\in \Z^2$
and all $\xi \in \{0,1\}^{\Z^2 \setminus \{x\}}$ we have
\[
\eps \, < \, \P(X(x)=1 \, | \, X(\Z^2 \setminus \{x\}) = \xi)
\, < \, 1-\eps \, .
\]
\end{defn}
In this paper, we show by means of concrete examples that
translation invariance together with finite energy is {\em not}
sufficient to rule out coexistence of an infinite
open cluster, and an infinite closed cluster. What's more, it does not even
help if we replace finite energy by
uniform finite energy:
\begin{thm}  \label{thm:main}
There exists a translation invariant site percolation on $\Z^2$ that
satisfies uniform finite energy and that produces a.s.\ an infinite open
cluster and an infinite closed cluster.
\end{thm}
It is a classical result of Burton and Keane \cite{BK} that
translation invariance and finite energy together are enough to rule out
the existence of more than one infinite open cluster (and, by
symmetry, more than one infinite closed cluster), so any example
witnessing Theorem \ref{thm:main} must have a.s.\ exactly one
infinite cluster of each kind.

The rest of this paper is devoted to examples exhibiting such coexistence.
The example witnessing Theorem \ref{thm:main} requires a somewhat
elaborate construction, and is therefore
postponed to Section \ref{sect:second}.
Along the way, we answer affirmatively 
(in Theorem \ref{thm:two_big_clusters}) the question of whether,
still assuming translation invariance, coexistence is possible
between an infinite open cluster and an infinite closed cluster that are
both robust under i.i.d.\ thinning. 
Before that, and in order to offer the reader
some intuition for the problem, we first present
a slightly less involved construction in Section \ref{sect:first},
which satisfies finite energy but not uniform finite energy.

\section{First construction} \label{sect:first}

The purpose of this section is to give an example which proves
the following weaker version of Theorem \ref{thm:main}.
\begin{prop} \label{prop:finite_energy}
There exists a translation invariant site percolation on $\Z^2$ that
satisfies finite energy and that produces a.s.\ an infinite open
cluster and an infinite closed cluster.
\end{prop}
The construction will be based on the notion of a uniform spanning tree
for the $\Z^2$ lattice, first studied by Pemantle \cite{Pem} and
later by Benjamini et al.\ \cite{BLPS} and others.

A {\bf spanning tree} of a connected graph $G=(V,E)$ is a
connected subgraph of $G$ that contains all vertices $v\in V$ but no
cycles. Any finite
such $G$ has a finite number of possible spanning trees, and
a {\bf uniform spanning tree} for $G$ is therefore elementary to define
in this finite setting: it is the random spanning tree for $G$ obtained
by choosing one of the possible spanning trees
at random according to uniform distribution. This procedure
may be identified with a probability measure $\mu$ on $\{0,1\}^E$ which we
call the uniform spanning tree measure for $G$.

If we now move on to the case where $G=(V,E)$ is infinite but locally finite,
the concept of a uniform spanning tree is less elementary, because there
may be infinitely many (even uncountably many) possible
spanning trees. Pemantle \cite{Pem} showed that the following natural
definition makes sense. By an {\bf exhaustion} of $G$, we mean a sequence
$G_1=(V_1, E_1), G_2=(V_2,E_2), \ldots$ of connected finite subgraphs
that exhausts $G$ in the sense that every $v\in V$ and every $e\in E$ is
contained in all but at most finitely many $G_i$'s. For each $G_i$ we know
how to pick a uniform spanning tree, so we may define a probability measure
$\mu_i$ on $\{0,1\}^E$ whose projection on $\{0,1\}^{E_i}$ is the uniform
spanning tree measure for $G_i$ while the projection on
$\{0,1\}^{E \setminus E_i}$ may be defined arbitrarily. It turns out
(see \cite{Pem}) that the $\mu_i$'s converge (in the product topology)
to a limiting measure $\mu$ on $\{0,1\}^E$.
Furthermore $\mu$ is
concentrated on subgraphs of $G$ consisting of a union of finitely
or infinitely many
infinite trees (i.e., not necessarily a single tree as might be tempting
to believe). Pemantle considered the case where $G$ is the $\Z^d$ lattice
-- having vertex set $V=\Z^d$ and edge set $E$ consisting of edges connecting
$L_1$-nearest neighbors --
and showed that the number of trees is a $\mu$-a.s.\ constant, equalling
$1$ for $d\leq 4$, and $\infty$ for $d\geq 5$. The case which concerns
us is $d=2$, where the resulting spanning tree, other than being unique,
also has the following interesting properties:
\begin{description}
\item{\bf One end.}
For every vertex $x \in \Z^2$, there exists $\mu$-a.s.\ exactly one infinite
self-avoiding path in the tree starting at $x$.
\item{\bf Self-duality.}
Consider the dual lattice $\tilde{\Z}^2$, with vertex set
$\tilde{V}=\Z^2=\Z^2+ (\frac{1}{2}, \frac{1}{2})$ and edge set $\tilde{E}$
consisting of edges connecting $L_1$-nearest neighbors. In the
natural planar embeddings of $(V,E)$ and $(\tilde{V}, \tilde{E})$,
each edge $e \in E$ crosses exactly one edge $\tilde{e} \in \tilde{E}$.
Suppose we pick $Y\in \{0,1\}^E$ according to $\mu$, and then pick
$\tilde{Y}\in \{0,1\}^{\tilde{E}}$ by declaring each $\tilde{e}\in \tilde{E}$
present in $\tilde{Y}$ if and only if the edge $e \in E$ that it crosses is
absent in $Y$. Then, it turns out, the distribution of $\tilde{Y}$ is
the same as that of $Y$ (apart from the $(\frac{1}{2}, \frac{1}{2})$
shift). In particular, $\tilde{Y}$ consists of a single one-ended
spanning tree for $\tilde{G}$.
\end{description}
Using Pemantle's spanning tree construction $Y$, we construct a site
percolation $X\in \{0,1\}^{\Z^2}$ as follows; it should be viewed as a picture
of $Y$ and $\tilde{Y}$ scaled up by factor $2$.
Writing $x\in \Z^2$ in terms of its coordinates as $x=(x_1,x_2)$
\[
X(x_1, x_2) = \left\{
\begin{array}{ll}
1 & \mbox{ if $x_1$ and $x_2$ are both even \, (these sites represent the
vertices of $Y$)} \\
0 & \mbox{ if $x_1$ and $x_2$ are both odd \, (these sites represent the
vertices of $\tilde{Y}$).} \\
\end{array} \right.
\]
The remaining sites respresent crossing pairs of edges in $E$ and $\tilde{E}$:
a $1$ indicates the presence of $e\in E$ and a $0$ that of its dual
edge $\tilde{e}$. More precisely,
if $x_1$ is even and $x_2$ is odd, we set $X(x_1, x_2)=1$ iff
the edge $e \in E$ linking $(\frac{x_1}{2}, \frac{x_2-1}{2})$ to
$(\frac{x_1}{2}, \frac{x_2+1}{2})$ is present in $Y$;
while if $x_1$ is odd and $x_2$ is even, we set $X(x_1, x_2)=1$ iff
the edge $e \in E$ linking $(\frac{x_1-1}{2}, \frac{x_2}{2})$ to
$(\frac{x_1+1}{2}, \frac{x_2}{2})$ is present in $Y$.

This defines $X \in \{0,1\}^{\Z^2}$. It is clear from the construction that
$X$ produces a.s.\ a single infinite open cluster, a single infinite closed
cluster, and no finite clusters. To serve as a couterexample proving
Proposition \ref{prop:finite_energy}, it is however deficient
in two ways, as (i) it fails to be translation invariant, and
(ii) it fails to exhibit finite energy.
Translation invariance is fixed by letting $\hat{X} \in \{0,1\}^{Z^2}$ equal
$X$ shifted by a random amount equalling $(0,0)$, $(0,1)$, $(1,0)$ or
$(1,1)$, each with probability $\frac{1}{4}$; the resulting
site percolation $\hat{X}$ is easily seen to be translation invariant.

To modify it again to
give it finite energy, note first that for each $x\in \Z^2$ which is
open in $\hat{X}$ there is a unique infinite self-avoiding
open path in $\Z^2$ starting at $x$ and using only open vertices,
and analogously for each closed $x\in \Z^2$. Thus, each $x\in \Z^2$
has the following property: removing $x$ would cut off a finite (possibly
$0$) number of vertices from either the infinite open or the infinite closed
cluster of $\hat{X}$. Write $b(x)$ for this number, and construct
$\bar{X}\in \{0,1\}^{\Z^2}$ by letting
\begin{equation}  \label{eq:perturbation}
\bar{X} = \left\{
\begin{array}{ll}
\hat{X}(x) & \mbox{with probability } 1-2^{-b(x)-1} \\
1- \hat{X} & \mbox{with probability } 2^{-b(x)-1}
\end{array} \right.
\end{equation}
independently for each $x \in \Z^2$. This defines $\bar{X}$, which
clearly satisfies finite energy. Furthermore, if $x$ is a site which
is open in $\hat{X}$, then the expected number of sites in the unique
infinite self-avoiding path from $x$ in $\hat{X}$ that flip in the
mapping (\ref{eq:perturbation}) is bounded by $\sum_{i=1}^\infty 2^{-i}$
and therefore finite. So we have a.s.\ that for some vertex in the path
and onwards, no vertex is flipped. Hence $\bar{X}$ has
an infinite open cluster. Similarly we get that it has an infinite closed
cluster. Thus, it has all the properties needed to warrant the
statement that {\bf Proposition \ref{prop:finite_energy} is established.}

\section{Second construction}  \label{sect:second}

Most of the work needed to prove Theorem
\ref{thm:main} is contained in the proof of the following Theorem
\ref{thm:two_big_clusters}. We will need some additional standard
terminology. For an infinite but locally finite graph $G$, define
the site percolation critical value
$p_{c,site}(G)$ to be the infimum over all $p\in [0,1]$ such
that i.i.d.\ site percolation on $G$ with retention parameter $p$ produces
a.s.\ at least one infinite open cluster. Also, let $p_{c,bond}(G)$ be the
analogous critical value for i.i.d.\ bond percolation on $G$, i.e., for the
percolation process where it is the edges (rather than the vertices)
that are removed at random.
The following result is well known; see, e.g.,
\cite[Thm. 1.1]{H00}.
\begin{lem}  \label{lem:site_vs_bond}
For any graph $G$ of bounded degree, we have $p_{c,site}(G)<1$ if and only
if $p_{c,bond}(G)<1$.
\end{lem}
Given a site
percolation $\hat{X}$ on $\Z^2$, we write $G_{open}(\hat{X})$ for the (random)
graph whose vertex set consists of all $x \in \Z^2$ such that
$\hat{X}(x)=1$, and whose
edge set consists of all pairs of such vertices at $L_1$-distance
$1$ from each other. Analogously, $G_{closed}(\hat{X})$
has vertex set consisting
of all $x \in \Z^2$ such that $\hat{X}(x)=0$, and edge set consisting
of all pairs of such vertices at $L_1$-distance $1$ from each other.
\begin{thm}  \label{thm:two_big_clusters}
There exists a translation invariant site percolation $\hat{X}$
such that with probability $1$ we have both $p_{c,bond}(G_{open}(\hat{X}))<1$
and $p_{c,bond}(G_{closed}(\hat{X}))<1$.
\end{thm}
Before proving this result, which is our main task, we show how
it easily implies Theorem \ref{thm:main}.

\medskip\noindent
{\bf Proof of Theorem \ref{thm:main} from Theorem \ref{thm:two_big_clusters}.}
Let $\hat{X}$ be as in Theorem \ref{thm:two_big_clusters}. By Lemma
\ref{lem:site_vs_bond}, we then have a.s.\ that
$p_{c,site}(G_{open}(\hat{X}))<1$ and $p_{c,site}(G_{closed}(\hat{X}))<1$.
We can then find an $\eps \in (0, \frac{1}{2})$ such that
\begin{equation} \label{eq:uniform_eps}
\P(p_{c,site}(G_{open}(\hat{X}))<
1 - \eps\, , \, p_{c,site}(G_{closed}(\hat{X}))<1 - \eps)\, > \, 0 \, .
\end{equation}
In fact, we may without loss of generality assume that the event
in (\ref{eq:uniform_eps}) has probability $1$, because the event
is translation invariant so that conditioning on
it does not mess up translation invariance.

Now obtain another site percolation $\bar{X}$ from $\hat{X}$ by letting,
for each $x \in \Z^2$ independently,
\begin{equation}  \label{eq:uniform_perturbation}
\bar{X} = \left\{
\begin{array}{ll}
\hat{X}(x) & \mbox{with probability } 1-\eps \\
1- \hat{X} & \mbox{with probability } \eps
\end{array} \right.
\end{equation}
It is immediate that the translation invariance property of $\hat{X}$
is inherited by $\bar{X}$. Furthermore, the transformation
(\ref{eq:uniform_perturbation}) implies (regardless of the details of
$\hat{X}$) that
\[
\P(\bar{X}(x)=1 \, | \, \bar{X}(\Z^2 \setminus \{x\})
\in [\eps, 1- \eps] \mbox{ a.s.,}
\]
so $\bar{X}$ satisfies uniform finite energy.
Next, $p_{c,site}(G_{open}(\hat{X}))< 1 - \eps$ implies that the set of sites
in $G_{open}(\hat{X})$ that remain unflipped through the transformation
(\ref{eq:uniform_perturbation}) contains an infinite cluster; and analogously
for $G_{closed}(\hat{X})$.
In summary, $\bar{X}$ has all the properties needed
to warrant Theorem \ref{thm:main}.
$\Cox$

\medskip\noindent
It remains to prove Theorem \ref{thm:two_big_clusters}. It is instructive
to think about why the $\hat{X}$ from Section
\ref{sect:first} will not do. In that example, for
each open vertex $x \in \Z^2$, $G_{open}(\hat{X})$ contains
only a single infinite self-avoiding path starting at $x$.
Carrying out i.i.d.\ bond percolation with retention parameter $1- \eps$
on $G_{open}(\hat{X})$ will, regardless of how small $\eps>0$ is,
a.s.\ kill at least one edge on this path and thus cut off
$x$ from any infinite cluster. Thus, $p_{c,bond}(G_{open}(\hat{X}))=1$ (and,
analogously, $p_{c,bond}(G_{closed}(\hat{X}))=1$), so this choice of
$\hat{X}$ fails to be a witness to Theorem \ref{thm:two_big_clusters}.

What made the infinite clusters of
$\bar{X}$ of Section \ref{sect:first} survive was the
inhomogeneity of the retention probabilities, sufficiently
rapidly approaching $1$
as we moved from $x \in \Z^2$ off along its single self-avoiding path
to infinity. When the retention parameter is set fixed at $1- \eps$,
we could try another approach: to replace the single path from
$x$ to infinity by a road that becomes
progressively broader (and therefore more robust to random thinning)
as we move along it.
Some intuitive evidence that this should be doable comes from the work
of Grimmett \cite{G81} and others concerning i.i.d.\
bond percolation on graphs $G_f$ arising
by restricting the $\Z^2$ lattice to vertices $x=(x_1, x_2)\in \Z^2$ with
$x_1\geq 0$ and $0 \leq x_2 \leq f(x_1)$
(and the usual nearest-neighbor edges connecting them), where
$f: \Z_+ \rightarrow \R_+$ is a function
that grows towards infinity as its argument goes to infinity.
It turns out that a relatively slow growth of $f$ suffices to
ensure that $p_{c, bond}(G_f)<1$; in particular, Grimmett showed that
the critical value $p_{c, bond}(G_f)$
equals that of the full $\Z^2$ lattice
(i.e., $p_{c, bond}(G_f)=1/2$) if and only if
$\lim_{n \rightarrow \infty} f(n)/\log(n)=\infty$.

The fact that such slow growth of $f$ is enough suggests that
it should be possible to modify the tree-structure of the
$\hat{X}$ of Section \ref{sect:first} in such a way as to obtain
a witness to Theorem \ref{thm:two_big_clusters}. This is what we set
out to do in the following. For technical reasons, we opt for a
tree-like structure with a lot more regularity than the example
in Section \ref{sect:first}. Our construction will be built up from
rectangular sets on a sequence of larger and larger scales. The
percolation theory developed in the last few decades offers an abundance of
results concerning crossing probabilities in i.i.d.\ percolation
on such rectangles. We will settle for one which is due to Bollob\'as and
Riordan -- see Lemma \ref{lem:BR} below -- although other choices
would certainly have been possible.

Due to the amount of work needed to prove Theorem \ref{thm:two_big_clusters},
we divide it into a number of smaller portions.
First, in Section \ref{sect:rectangles}, we
introduce the terminology needed for a precise discussion of
crossing probabilities for i.i.d.\ percolation and the Bollob\'as--Riordan
result. Then, in Section \ref{sect:building_blocks}, we define the basic
building blocks of our construction. In Sections \ref{sect:open} and
\ref{sect:closed}, we go on to some preliminary considerations that will
be crucial for showing $p_{c,bond}(G_{open})<1$ and $p_{c,bond}(G_{closed})<1$,
respectively. The construction is completed in
Section \ref{sect:actual_construction}, and in
Section \ref{sect:finishing_the_proof} we finally establish
$p_{c,bond}(G_{open})<1$ and $p_{c,bond}(G_{closed})<1$,
thus completing the proof.

\subsection{Rectangles and crossing probabilities}  \label{sect:rectangles}

{{{If $S_1, S_2$ are subsets of $\Z^{d}$ then we will call them {\bf congruent}
if there exists a $v \in \Z^{d}$ so that $S_1= S_{2}+v$. Note that for us $d$ will be $1$ or $2$.}}}
If $B \subset {\Z} $ and $B= (a,b)  \cap {\Z} $, then we say that
 $B$
is a {\bf block}. A subset $R$ of $\Z^{2}$ will be called a {\bf
rectangle} if it can be written as $ R = B_{1} \times B_{2} $, where
the $B_i$'s are blocks; if the blocks are congruent, we call it a
{\bf square}. If $ B_1$ has $l$ elements and $B_2$ has $k$, we say
that $R$ is an $l \times k$ rectangle. The sets $ ( \min B_{1})
\times B_{2}$ and $ ( \max B_{1}) \times B_2 $  are called the
({\bf left} and {\bf right}, respectively)
{\bf vertical sides}  of $R$. The sets
$B_{1} \times ( \min B_2) $ and $ B_{1} \times (\max B_2) $  will be
the ({\bf bottom} and {\bf top}, respectively)  {\bf horizontal sides}.

We shall need the notion of
{\bf crossing} in a rectangle when preforming
i.i.d.\ bond percolation with retention probability $p$
-- indicated by writing $\P_p$ for the probability measure --
on it. For such a percolation process on a rectangle $R$, the event
$H(R)$ defined as the set of those  subgraphs of $R$ containing a
path between the two different vertical sides will be called a {\bf
horizontal crossing} in $R$. The event $V(R)$ which we define by
interchanging the words vertical and horizontal above called a {\bf
vertical crossing } in $R$.

Furthermore, we say that a rectangle $Q = A_{1}
\times A_2$ is {\bf well-joined} to the rectangle $R=B_{1} \times
B_2$ if {\em either}
\begin{quote}
$ A_{1} \subset B_{1}$ and $B_{2} \subset A_2 $
(in which case we say
that their being ``well-joined" is of type {\bf vertical to
horizontal} or $V \rightarrow H$)
\end{quote}
{\em or}
\begin{quote}
$ B_{1} \subset A_{1}$ and $A_{2} \subset B_2$
(in which  case we say that their being  ``well-joined" is of type {
\bf horizontal to vertical} or $H \rightarrow V$).
\end{quote}
 If $R_{1} , R_{2} ,
\dots , R_{m}, \dots $ is a sequence of rectangles, we say that it
is well-joined if every pair of consecutive rectangles from the
sequence is well-joined and the sequence of their types is
alternating (i.e: $ \dots V \rightarrow H, H \rightarrow V, V
\rightarrow H, H \rightarrow V, \dots $ ).  The importance
of this concept will be the following: if we have a sequence of
well-joined rectangles  $R_{1} , R_{2} , \dots , R_{m}, \dots $ and the
first type is (say) $V \rightarrow H$  and if all the events
$$V(R_{1}), H(R_{2}), V(R_{3}), H(R_{4}), \dots ,V(R_{2k-1}),
H(R_{2k}), \dots
$$ hold then we can easily extract an infinite path from the
individual crossings (given the appropriate vertical or horizontal
crossings for the rectangles in the sequence). Moreover, if we know a
lower bound for the individual probabilities of the above events,
then by the well-known Harris--FKG inequality which states
that for i.i.d.\ percolation any two increasing events are
positively correlated (see, e.g.,
\cite[Thm.\ 2.4]{G}) we can get a lower bound for the
probability of an infinite path simply by multiplication. We shall
make use of the following  result of Bollob\'as and Riordan \cite{BR}.
\begin{lem}  \label{lem:BR}
Fix an integer $\lambda  > 1 $ and a  $p\in (\frac{1}{2},1)$.
We can then find constants $\gamma= \gamma(\lambda, p) > 0$ and $n_0
=n_{0}(\lambda,p)$ such that if $n> n_0$,
then for each $ \lambda n \times n$ rectangle $R$ we have
$\P_{p} (H (R))> 1 - n^{- \gamma } $.
\end{lem}
From this result, we obtain the following.
\begin{cor}  \label{cor:BR}
For any fixed $p\in (\frac{1}{2},1)$, we can find
a $c>0$, a positive integer $n_0$ and a  $ \gamma
>0 $ such that for any positive integer $L$ the following holds.
If $R$ is an $Ln \times n$ rectangle where $ n > n_0$,
then $\P_{p}(H(R)) \geq c^{L/{n^{\gamma}}} $.
\end{cor}

\medskip\noindent
{\bf Proof.} Take $ \lambda =3$ in Lemma \ref{lem:BR} and let
$\gamma$ be the corresponding $\gamma(3,p) $ and $n_0$ be the
corresponding $n_0 (3,p)$. Let $R$ be an $Ln \times n$ rectangle.
$R$ can be covered by overlapping ``little" $3n \times n$ rectangles
in such a way that the intersection of a consecutive pair of them is
an $n \times n$ square and we can do it in such a way that
altogether the number of the $3n \times n$ {{{rectangles}}} and $n \times n$
{{{squares}}} is not greater than $L$. Notice that if we have
horizontal crossings for all the $3n \times n$ rectangles and
vertical crossings for all the $n \times n$ squares, then we have a
horizontal crossing for the whole $Ln \times n$ rectangle. Then, by
the Harris--FKG inequality, we can estimate $\P_{p} (H(R))$ from
below as $(1-n^{-\gamma})^L$. But this quantity equals
$((1-n^{-\gamma})^{n^{\gamma}})^{L/{n^{\gamma}}}$, so the corollary follows
from the fact that $(1-n^{-\gamma })^{n^{\gamma}}$ is bounded away
from zero. $\Cox$

\subsection{Building blocks}  \label{sect:building_blocks}

For a finite set $ K \subset \Z $, we let $ {\bf conv}(K)$ denote the smallest
block containing $K$. If  ${{{C}}}$ has the form
$${{{C}}}= \bigcup_{k \in {\Z} } (B + (l+d)k) $$  where $d>1$ is
some integer, and $B$ is a block with $ |B|=l
$, then we say that ${{{C}}}$ is a {\bf block progression}, and we refer
to $l$ as the {\bf block length} and $d$ as the {\bf block distance}
in ${{{C}}}$. We say that $(l, d)$ is the {\bf parameter } of ${{{C}}}$.  We
will refer to the sets $B_{k}=B+(l+d)k$ as the {\bf blocks} of ${{{C}}}$.
We call a block $D$ a {\bf gap} of ${{{C}}}$ if it is in the complement of
${{{C}}}$ and maximal with that property. Note that in that case $|D|=d$.
If $L$ is a positive integer, $ T\subset \Z $ and ${{{C}}}$ is a block
progression as above, then we say that $T$ is a { \bf block
progression over} ${{{C}}}$ {\bf with factor} $L$ if
$$T= \bigcup _{k \in \Z} ( D+(l+d)Lk) $$ where $D$ is a gap of  ${{{C}}}$.
Let $C_1$ and $C_2$ be two congruent block progressions. The blocks
of $C_i$ will be denoted as $B_{j}^{i}$ where $ j \in \Z$.
 Let
$$V_j= B_{j}^{1} \times \Z \hbox{ and } H_j= \Z
\times B_{j}^{2} \, . $$
 Then the set $G \subset {\Z}^2$ defined as
 $$ G= {\left( \bigcup_{k \in \Z} V_k \right) } \cup { \left(
\bigcup_{j \in \Z} H_j \right) }$$
 will be called the {\bf grid determined by} $C_1$ and $C_2$.
 The {\bf parameter of the grid} above will be the parameter of
 $C_i$.
 If $G$ and $H$ are grids we say that
 $H$ is a {\bf grid over} $G$ with factor $L$ if,
 whenever $G$ is determined by $C_i$ and $H$ is determined by $T_i$
 for $i \in \{ 1,2 \} $, $T_i$ is a block progression over
$C_i$ with factor $L$.

We now go on to define finite analogues of the above concepts.
If $B$ is a block and $C=  \bigcup_{k=0}^{q-1} (B+(l+d)k) $, then we
say that $C$ is a {\bf block complex}.  Next we define the notion of
a {window}. Let $C_1$ and $C_2$ be two congruent block complexes.
The blocks of $C_i$ will be denoted as $B_{j}^{i}$, where $ j \in
\{0,1, \dots , {q-1} \} $. Let $$V_j= B_{j}^{1} \times {\bf
conv}(C_2) \, \hbox{ and } \, H_j= {\bf conv}(C_1) \times B_{j}^{2} \, .$$
Then the set $W \subset {\Z}^2$ defined as $$ W= {\left(
\bigcup_{k=0}^{q-1} V_k \right) } \cup {\left( \bigcup_{j=0}^{q-1} H_j \right) }$$
will be called the {\bf window determined by} $C_1$ and $C_2$.

We shall call the $V_i$'s and $H_j$'s the {\bf frames} of the given
window. The convex hull of a window in $ \R^2 $ is a square whose
intersection with $ \Z^2 $ is the {\bf shade } of the window. If we
take the set theoretic complement of the frames in the shade, then
the resulting set  splits into squares in $ \Z^2 $. We refer to
those squares as the {\bf panes} of the window.
For a window $W$ as above let us refer to the corresponding block
length (independent of $i$ and $j$)  $|B_{j}^i|$ as the ``frame
width" of $W$, denoted as ${\bf w}(W)$. Also $|{\bf conv}(C_{i})|$
will be called the ``side length" of W and we denote it as ${\bf
s}(W)$.

 For a window $W$ as above we define its {\bf fork} as
follows: It will be the union of $q-1$ vertical parts and one
horizontal part. The vertical  parts (we shall call them the {\bf
cut-frames} of the fork) are the sets of the form
 $$\cut(V_i) := V_i \setminus  H_{q-1} $$ where  $ i \in \{1, \dots ,q-1 \}$.
 That is,  we cut off each vertical strip at the top and we throw away
the leftmost
vertical strip. The horizontal part (which we shall call the {\bf bottom} of
the fork) will be
$$H_{0} \setminus V_0 \, . $$
Thus altogether the fork of $W$ is defined as
$$ F(W):= \left(\left(  \bigcup_{i=1}^{q-1}
V_i \right)  \cup H_{0} \right) \setminus
\left( V_0 \cup H_{q-1} \right) \, .$$

\subsection{Preliminaries for $p_{c,bond}(G_{open})<1$} \label{sect:open}

 If
we have two windows $W$ and $W^{+}$, we write $W \prec W^{+}$ to
indicate that the shade of $W$ is a pane of $W^{+}$ (note that this
relation is not transitive). If we have a sequence $$ {\cal S}= W_1,
\dots ,W_k, \dots
$$ of windows then we write $W_{1} \prec W_{2} \prec \dots \prec
W_{k} \prec W_{k+1} \prec \dots $ to indicate $W_k \prec W_{k+1}$
for each $k$.

If we have a sequence as above we define the set $\ERB_{k}( \cal
S)$, and when the sequence $ \cal S$ is understood  we write
simply $\ERB_{k}$; $\ERB$ stands for ``Escape route to
the Right and to the Bottom''.
Note that for $k>1$, $\ERB_{k}$ will be the union of two
rectangles.
 First observe that if $W \prec W^{+}$ holds, then there is a unique
 vertical frame $V^{+}$ of $W^+$ which is attached to $W$ from the
right in the
sense that
  $(W+(1,0)) \cap V^+$ is
 nonempty.
 Now consider $W_{k-1} \prec W_{k} \prec W_{k+1}$. Let $V_{k}^+$ be the unique
 vertical frame attached to $W_{k-1}$ from the right. Then
 $\cut(V_{k}^+)$ will be one of the rectangles whose union is
 $\ERB_{k}$. To define the other rectangle we take the bottom
 $B_k$
 of $F(W_k)$ and we extend it to the right to get the ``extended
 bottom''
 $$E_k := \bigcup_{j=0}^{{\bf w}(W_{k+1})}(B_{k}+(j,0))\, .$$
 Now let us define
\[
\ERB_{k}:=\cut(V_{k}^+) \cup E_k \, .
\]
 Note that
\begin{equation} \label{eq:star}
\cut(V_{k}^+) \mbox{ is a } {\bf {w}}(W_k) \times ({\bf s}(W_k)-{\bf
 {w}}(W_k)) \mbox{ rectangle}
\end{equation}
while
\begin{equation}  \label{eq:star_star}
E_k \mbox{ is a }
({\bf s}(W_k)-{\bf {w}}(W_k)+{\bf {w}}(W_{k+1})) \times {\bf {w}}(W_k)
\mbox{ rectangle.}
\end{equation}
Now we want to extend this definition to $k=1$ as well.
 Note that the definition for $E_k$ can be adapted to the case $k=1$
 with no difficulty. The only thing that we do not have a natural
 choice for is a cut-frame. We simply define $\ERB_{W_1}$ as $F(W_1)
 \cup E_1$.
 Finally we define the {\bf road} ${\bf r}({\cal S})$
of the sequence $\cal S$ above as
 $${\bf r}({\cal S}) := {\bigcup_{k=1}^{\infty}} \ERB_{k} \, .$$
 The importance of the road is the following.
 If in the i.i.d.\
percolation each edge inside $\ERB_{1}$  remains open and for each
$\ERB_{k}$ for $k>1$ we have a vertical crossing for the
corresponding cut-frame and a horizontal crossing for the
corresponding bottom, then for each point of $F(W_1)$ there is an
open path to infinity. Note that, besides the exceptional $k=1$ case,
the remaining parts of the road can be considered as a sequence of
well-joined rectangles.

\subsection{Preliminaries for $p_{c,bond}(G_{closed})<1$} \label{sect:closed}

If we are in the shade of a window but
not in its fork, then we can move to the left top corner of the
window by moving always outside of the fork.
More specifically, if we have
$$ {\cal S} =W_{1} \prec W_{2} \prec \dots \prec
W_{k} \prec W_{k+1} \prec \dots $$ then we define the corresponding
$\ELT_{k}(\cal S)$ as follows, $\ELT$ being short for
``Escape route to the Left and to the Top''. Consider the pair $W_{k-1}
\prec W_k$. Take the leftmost vertical frame  $V_{k-1}^l$ of
$W_{k-1}$ note that this is contained in the complement of
$F(W_{k-1})$. Let $ \ext(V_{k-1}^l)$ be the rectangle maximal for
the following properties.
 It is contained in the shade of $W_k$, while {{{its ``horizontal component" is the same as that of $V_{k-1}^l$ in the sense that if $V_{k-1}^l=A \times B$ and $ \ext(V_{k-1}^l)={\hat A} \times {\hat B}$ then $A={\hat A}$, and we also have}}}
$$V_{k-1}^l \subseteq \ext(V_{k-1}^l)$$
and
$$ \ext(V_{k-1}^l) \subseteq {(F(W_k) \cup F(W_{k-1}))}^C \, .$$
Also let $H_k^t$ be the topmost horizontal frame of $W_k$ and let
$$\ELT_k := \ext(V_{k-1}^l) \cup H_k^t \, . $$
For the record, note the size of these two rectangles:
$\ext(V_{k-1}^l)$ is a ${\bf w}(W_{k-1}) \times ({\bf s}(W_k)-{\bf
w}(W_k))$ one, while $H_k^t$ is a ${\bf s}(W_k) \times {\bf w} (W_k)$
one.
 If we take a similar union for the $\ELT$'s as we had for the
$\ERB$'s then we will have an infinite ``road'' to infinity moving
strictly outside of the forks of the windows in the sequence (but
still in the windows).

\subsection{The actual construction}  \label{sect:actual_construction}

The site percolation $\hat{X} \in \{0,1\}^{\Z^2}$ that we are about to
define will depend on two initial parameters $d_{0}$
and $l_{0}$ and a sequence of positive integers $ L_{1}, \dots
,L_{k}, \dots $ where the latter sequence ``grows fast" in a
later-specified way. We choose $d_0 \geq l_0 > n_0$ where $n_0$ is
from Corollary \ref{cor:BR}.

Note that there are only finitely many different translates of a
given grid so  we can choose uniformly a grid $G_0$ with parameter
$(l_{0} , d_{0})$ among the finitely many congruent copies.  If
$G_k$ has been defined for  a positive integer $k$, then let
$G_{k+1}$ be a uniformly chosen grid over $G_k$ with factor
$L_{k+1}$. If the $L_i$ grow fast enough, then a.s. any $x \in \Z^2$
will be in $G_k$ for only finitely many $k$.

If the grid
$G_i$ has parameter $(l_i , d_i)$, then $G_{i+1}$ will have
parameter $(l_{i+1},d_{i+1})=(d_{i}, L_{i+1}l_i+(L_{i+1}-1)d_{i})$.
Then for $l_{i}+d_{i}$ we have the simple recursion
$$ l_{i+1}+d_{i+1} = L_{i+1} (l_i+d_i)$$
which clearly implies
\begin{equation}  \label{eq:main_recussion}
l_{i+1}+d_{i+1} = ( \prod_{j=1}^{i+1}L_{j})(l_0+d_0) \, .
\end{equation}
Now color the points $x$ of $\Z^2$ with colors $ -1, 0, 1, 2, \dots
,n, \dots $ as follows: if $x$ is not in any of the $G_{k}$, then
$x$ gets the color $-1$, otherwise it gets the largest $k$ for which
$x \in G_{k}$.

It is crucial to make sure that any vertex be in only
finitely many of the $G$'s, for which a Borel--Cantelli argument is
enough if the $L_k$'s grows fast enough. We now give a sufficient
condition for that. Let us estimate the probability that the origin
is in $G_k$ (by invariance the same estimate works for any given
vertex). If we have a grid $H$ with parameter $(l,d)$ then instead
of looking at this as a union of certain vertical and horizontal
``infinite rectangles" we can visualize $\Z^2$ as {{{partitioned into a}}} disjoint union
of $(l+d) \times (l+d)$ squares {{{and consider the portion $H$
has within each of the squares. These portions will give us the probabilities that a particular point is contained in $H$. To compute these portions
 we choose the squares so that their intersection with $H$
 is especially simple, namely for each square $K$ from the partition the following holds:
 $H\cap K$ is the union of two rectangles $R_v$ and $R_h$ (here $h,v$ refers to ``horizontal" and ``vertical" respectively) so that $R_v$ has type $l \times (l+d)$ and $R_h$ has type
 $(l+d) \times l$ and $R_v$ is the ``leftmost" rectangle of that type contained in $K$ while $R_h$ is the ``topmost" one, meaning that neither $(-1,0) +R_v$ nor $R_h+(0,1)$ is contained in $K$. This gives us that $|H\cap K|=ld+l(l+d)$ while obviously $K=(l+d)^2$.}}}
  Then the probability of
the origin being in $H$ (if $H$ is uniformly selected as was the
case with the $G_k$'s) equals
$$(ld+l(l+d))/{(l+d)^2} \, .$$
Now let
us check what condition on $L_1, L_2, \dots $ needed to make the
Borel--Cantelli argument work. In order to do that consider
$H=G_{i+1}$ so the probability of the origin being in $H$ is
$$(l_{i+1}d_{i+1}+l_{i+1}(l_{i+1}+d_{i+1})) / {(l_{i+1}+d_{i+1})^2} \, ,$$
which, with a little bit of arithmetic, becomes
$$2l_{i+1}/(l_{i+1}+d_{i+1})- {l_{i+1}^2}/{(l_{i+1}+d_{i+1})^2} \, .$$
For Borel--Cantelli to work we need that summing these positive
numbers over $i$ yields a finite value, and for that it is clearly
enough that the sum of $l_{i+1}/(l_{i+1}+d_{i+1})$ converges. To see
how it relates to $L_1, L_2, \dots $ we spell out our recursions
again:
$$l_{i+1}/(l_{i+1}+d_{i+1})=d_{i}/(L_{i+1}(l_{i}+d_{i}))<1/{L_{i+1}} \, .$$
So it is enough to have $$\sum_{i=1}^{\infty } 1/{L_{i}} <\infty \, .$$

After this Borel--Cantelli interlude, we now turn back to the construction.
Observe that each color class splits into a disjoint union of
windows. Actually a more precise { \bf  ``structural observation"}
is true:
 A point $x$  of color class $k$ is always contained in a window $ W(x)$
 each  of whose points has  the same color with
 ${\bf w}(W(x))=l_k$ and ${\bf s}(W(x))=d_{k+1}$.
 Also for this $W(x)$ there exists a $W^{+}(x)$ each of whose
 points has color $k+1$ so that $W(x) \prec W^{+}(x)$.
 Altogether we find that for an $x$ of color class $k$ we have a
 sequence
 $W_{1}(x)
\prec W_{2}(x) \prec \dots \prec W_{j}(x) \prec W_{j+1}(x)  \prec
\dots $ of windows where each point of $W_j(x)$ is of color class
$k+j-1$.

 The construction is simply to take the forks of all of  the
windows: our translation invariant site percolation $\hat{X}\in \{0,1\}^{\Z^2}$
arises by assigning value $1$ to $x$ precisely for those $x \in \Z^2$
that belong to such a fork. Let us write $\RF$ (short for ``Random Forks'')
for $G_{open}(\hat{X})$ and $\RF^*$ for $G_{closed}(\hat{X})$.

\subsection{Nontriviality of the critical values}
\label{sect:finishing_the_proof}

It remains to show that $p_{c,bond}(\RF)<1$ and $p_{c,bond}(\RF^*)<1$;
we begin with the former. If $x \in \Z^2$ is in
$\RF$ consider  $W_{1}(x)
\prec W_{2}(x) \prec \dots \prec W_{j}(x) \prec W_{j+1}(x)  \prec
\dots $ as above.
  With positive probability $x$ is of color class
 $1$. Moreover, still with positive probability, each edge in $\ERB_1$
remains open.
 We will condition on this event.

Then we can use the notion of road ${\bf r}(x)$ introduced in
Section \ref{sect:open}.
Because of the conditioning we just declared,
we can focus on the parts in the road which
corresponded to indices $k>1$. Let us apply the strategy we
described in Section \ref{sect:rectangles} in connection
with the notion of being well-joined. We
need the side lengths of the rectangles constituting the road,
and to substitute the frame widths and side lengths of $\RF$
into the formulas (\ref{eq:star}) and (\ref{eq:star_star}) in Section
\ref{sect:open}.
The vertical rectangle at the $i$th step of the
road for $x$ is an $l_i \times (d_{i+1} -l_i)$ one while the next
horizontal one is a $(d_{i+1}+l_{i+1}-l_{i}) \times l_i$ one. Note
that from the point of view of crossing an ${2d_{i+1}} \times l_i$
(horizontal) rectangle and an $l_i \times 2d_{i+1}$  (vertical) is
just worse than any of the above so if we find a lower bound for
their having the appropriate crossings then that bound works for the
original rectangles as well.

Now by using the recursion (\ref{eq:main_recussion})
we obtain estimates
for the side lengths: $$ 2d_{i+1} > d_{i+1}+l_{i+1} = \left(
\prod_{j=1}^{i+1}L_{j} \right)(l_0+d_0)
>d_{i+1} \, .$$ Note that (simply because $l_{i+2}=d_{i+1} $) we also
have $$ 2l_{i+2} > \left( \prod_{j=1}^{i+1}L_{j} \right)(l_0+d_0) \, ,$$
and furthermore $$
2L_{i+1}L_{i}l_{i}=2L_{i+1}L_{i}d_{i-1}
> L_{i+1}L_{i}(d_{i-1}+l_{i-1})=d_{i+1}+l_{i+1}>d_{i+1} \, .$$ In other
words we have $2L_{i+1}L_{i} > d_{i+1}/{l_{i}}$.
 Now apply Corollary \ref{cor:BR} to
the above $l_i \times 2d_{i+1}$ rectangle $R$. Then $4L_{i+1}L_{i}$
may play the role of $L$ in the corollary, which then tells us that
\begin{equation}  \label{eq:each_level}
\P_p(V(R)) > c^{4L_{i+1}L_{i}/{l_{i}^{\gamma }}} \, .
\end{equation}
We next use the fact that the sequence of rectangles defined above (i.e.\
the ``road'' we get  when we take a vertical strip from the fork and
go down to the bottom horizontal one and the move to the vertical
strip in the next level and so on...) is well-joined. The estimate
(\ref{eq:each_level}) together with the Harris--FKG inequality implies
that the
probability of the sequence containing an infinite path is greater than
\begin{equation}  \label{eq:product}
\prod_{i=2}^{ \infty  }c^{2({4L_{i+1}L_{i}/{l_{i}^{\gamma }}})}
\end{equation}
where the factor $2$ in the power corresponds to taking both the
horizontal and vertical rectangles into account at a given step, and
the index $i$ going from $2$ corresponds to the conditioning
declared at the beginning of Section \ref{sect:finishing_the_proof}. The
product (\ref{eq:product}) is positive exactly when
\begin{equation}  \label{eq:what_we_want}
\sum_{i=2}^{\infty}L_{i+1}L_{i}/{l_{i}^{\gamma }} < \infty \, .
\end{equation}
Recall the balance we need to establish: on one hand,
the $L$'s need to grow
fast enough so that  Borel--Cantelli applies  to show only finitely many of
the events $x \in G_k$ hold, while
on the other hand they need to grow slowly enough to
make sure that the sum (\ref{eq:what_we_want})
converges. But of course with the given
conditions there is {\em plenty} of room for that because as we saw
the $l_i$ is essentially the product of all $L_k$'s up to index $i$.
We can even  allow the $L$'s to grow exponentially. Indeed, let $L_i=2^i$.
Then we see that the term corresponding to  index $i+2$ of the above
sum  is $ 2^{2i+5}/{l_{i+2}^{\gamma}}$. Now note that
$$ 2l_{i+2} > ( \prod_{j=1}^{i+1}L_{j})(l_0+d_0)=
2^{({(i+1)(i+2)}/2)}(l_0+d_0) \, .$$ We note that
$2^{2i+5}/{l_{i+2}^{\gamma }}$ can be bounded from above as some
constant multiplied by $ 2^{- \alpha {i^{2}} + \beta i + \delta }$
(where $\alpha>0, \beta, \delta \in \R $) whose sum (over $i$) is
clearly convergent. (In fact, we could consider even faster growing $L$'s
as long as we make sure that the product of the first some terms should be
much bigger than the next two terms.)

This justifies our claim that
$p_{c,bond}(\RF)<1$ for $\RF$, and it remains only to establish
the analogous claim $p_{c,bond}(\RF^*)<1$.
For that purpose we do a computation very similar to the above one
but now applied to the road defined by the $\ELT$'s. Note that the
sizes of the vertical and horizontal rectangles in $\ELT_k$ are $
d_{k+1} \times l_k $ for the horizontal one and $ l_{k-1} \times
(d_{k+1}-l_k) $ for the vertical one, and furthermore
that in this case both crossing probabilities for the above
considered two rectangles is not less then the horizontal crossing
probability for a $d_{k+1} \times l_{k-1} $ one.

First we need an estimate for the ratio $d_{i+1}/l_{i-1}$. We
use again the basic recursion for the $(l+d)$'s we had at the
``structural observation":
$$2L_{i+1}L_iL_{i-1}l_{i-1}=2L_{i+1}L_iL_{i-1}d_{i-2}>
L_{i+1}L_iL_{i-1}(d_{i-2}+l_{i-2})=d_{i+1}+l_{i+1}>d_{i+1} \, .$$
So now the quantity
$2L_{i+1}L_iL_{i-1}$ can play the role of $L$ from Corollary \ref{cor:BR}.
So we need
\begin{equation}  \label{eq:what_we_need_this_time}
\sum_{i=1}^{\infty}L_{i+2}L_{i+1}L_{i}/{l_{i}^{\gamma }} < \infty  \, .
\end{equation}
Now if we make the same kinds of estimates as for $\RF$,
we see that the $i$'th term in this case will be
$2^{3i+3}/{l_{i}^\gamma}$. So in the numerator we still
have an exponent linear in $i$, while in the denominator we have an
exponent of second order, so the sum in (\ref{eq:what_we_need_this_time})
is indeed finite,
and the proof of Theorem \ref{thm:two_big_clusters} is complete.

\medskip\noindent
{\bf Acknowledgement.} We thank Russ Lyons for valuable discussions.

\end{document}